\documentstyle[12pt]{article}
\newcommand{\Sub}{\mathop{\rm Sub}\limits}

\newcommand{\SSub}{\mathop{\rm SSub}\limits}

\newcommand{\const}{\mathop{\rm const}\limits}

\newcommand{\Var}{\mathop{\rm Var}\limits}

\newcommand{\supp}{\mathop{\rm supp}\limits}

\begin{document}

 \begin{center}

{\bf SUBGAUSSIAN AND STRICTLY SUBGAUSSIAN \\

\vspace{4mm}

RANDOM VARIABLES} \\

\vspace{4mm}

{\sc Eugene Ostrovsky, Leonid Sirota}\\

\vspace{3mm}

 Bar - Ilan University,  59200, Ramat Gan, ISRAEL; \\

 e-mail: eugostrovsky@list.ru \\
 e-mail: sirota3@bezeqint.net \\

\vspace{5mm}

\hspace{55mm} {\it Devoted to the memory of V.V.Buldygin}\\

 \vspace{5mm}

        {\bf Abstract}

\end{center}

\vspace{3mm}

  We study in this report the so-called Strictly Subgaussian (SSub) random variables (r.v.), which form  a very interest subclass
 of Subgaussian (Sub)  r.v., and obtain the exact exponential bounds for tail of distribution for sums of independent and disjoint such a variables,
 not  necessary to be identical distributed, and give some new examples  of SSub variables  to show  the exactness of our estimates. \par
  We extend also these results on the case of sums of subgaussian martingale differences, and show that the mixture of (Strictly)
  Subgaussian r.v. forms also (Strictly) subgaussian variable. \par

\vspace{3mm}


\vspace{3mm}

 \section { Introduction. Definitions. Notations. Examples.}

\vspace{3mm}

 Let $  \{\Omega, B, {\bf P}  \} $ be some non-trivial probability space with expectation $ {\bf E.} $ \\

 \vspace{2mm}

 {\bf Definition  1.1.}\\

 \vspace{2mm}

  We say that the {\it centered:} $ {\bf E} \xi = 0 $ numerical random variable (r.v.)
 $ \xi = \xi(\omega), \ \omega \in \Omega $ is subgaussian, or equally, belongs to the space $ \Sub(\Omega), $
 if there exists some non-negative constant $ \tau \ge 0 $ such that

$$
\forall \lambda \in R  \ \Rightarrow
{\bf E} \exp(\lambda \xi) \le \exp[ \lambda^2 \ \tau^2 /2]. \eqno(1.1).
$$

\vspace{3mm}

 The minimal value $ \tau $ satisfying (1.1) is called a  {\it subgaussian  norm}
of the variable $ \xi, $ write

 $$
 ||\xi||\Sub = \inf \{ \tau, \ \tau > 0: \ \forall \lambda \in R \ \Rightarrow {\bf E}\exp(\lambda \xi) \le \exp(\lambda^2 \ \tau^2/2) \}.
 $$

 Evidently,

$$
||\xi||\Sub = \sup_{\lambda \ne 0} \left[ \sqrt{ 2 \ln {\bf E}  \exp (  \lambda \xi)  }/|\lambda| \right].  \eqno(1.2)
$$

 This important notion was introduced by  J.P.Kahane \cite{Kahane1}; V.V.Buldygin and Yu.V.Kozachenko in \cite{Buldygin1} proved
that the set $ \Sub(\Omega) $  relative the norm $  ||\cdot|| $ is complete Banach space which is isomorphic to subspace
consisting only from the centered variables of Orlicz's space over $ (\Omega, B,P)  $ with $ N \ -$ Orlicz-Young function
 $ N(u) = \exp(u^2) - 1; $  see also \cite{Kozachenko1}.  \par

  For instance, let us consider the centered indicator (binary) random variable $ \xi_p: $

$$
{\bf P}(\xi_p = 1 - p) = p; \ {\bf P}(\xi_p = -p) = 1 - p, \ p \in (0,1);
$$
then (in our definitions and notations)

$$
||\xi_p||\Sub = \sqrt{ \frac{1 - 2p}{2 \ln(( 1 - p )/p)} }.
$$
 This important result was obtained independently in the articles \cite{Kearns1},  \cite{Buldygin3};  see also
 \cite{Berend1},  \cite{Schlemm1}, \cite{Ostrovsky2}.\par

\vspace{3mm}

  The detail investigation of this class or random variables with very interest  applications into the theory of random fields
  reader may found in the book \cite{Buldygin2};  we reproduce here some main facts from this monograph.\par

   If $ ||\xi||\Sub = \tau \in (0,\infty),  $ then

 $$
 \max [{\bf P}(\xi > x),  {\bf P}(\xi < -x)  ] \le \exp(- x^2/(2 \tau^2)  ), \ x \ge 0; \eqno(1.3)
 $$
 and  the last inequality is in general case non-improvable.  It is sufficient for this to consider the case when
 the r.v. $  \xi  $ has the centered Gaussian non-degenerate distribution.\par

  Conversely, if  $ {\bf E} \xi = 0 $ and if  for some positive finite constant $  K  $

 $$
 \max [{\bf P}(\xi > x),  {\bf P}(\xi < -x)  ] \le \exp(- x^2/K^2  ), \ x \ge 0,
 $$
 then $ \xi \in \Sub(\Omega) $ and $ ||\xi||\Sub < 4 K. $ \par

 The subgaussian norm in the subspace of the centered r.v. is equivalent to the following Grand Lebesgue Space (GLS)
 norm:

 $$
|||\xi||| := \sup_{s \ge 1} \left[ \frac{|\xi|_s}{\sqrt{s}} \right], \hspace{6mm} |\xi|_s =  \left[ {\bf E} |\xi|^s \right]^{1/s}.
 $$

 For the non-centered r.v. $ \xi $  the subgaussian norm may be defined as follows:

 $$
 ||\xi|| \Sub := \left[  \left\{ ||\xi - {\bf E} \xi||\Sub \right\}^2 + ( {\bf E} \xi)^2  \right]^{1/2}.
  $$

 More detail investigation of these spaces see in the monograph \cite{Ostrovsky1}, chapter 1.  \par

 Denote in the sequel for brevity  for any r.v.  $ \eta $

 $$
 \sigma^2(\eta) = \sigma^2 = \Var{\eta} = {\bf E}\eta^2 - ( {\bf E} \eta)^2.
 $$

\vspace{3mm}

{\bf Definition 1.2.} (See \cite{Buldygin2}, chapter 1.) \par

The subgaussian r.v. $ \xi $ is said to be {\it Strictly Subgaussian  (SSub), } iff

$$
\forall \lambda \in R \ \Rightarrow {\bf E} e^{\lambda \xi} \le e^{\lambda^2 \sigma^2(\xi)/2}, \eqno(1.4)
$$
or equally

$$
||\xi ||\Sub \le \sigma(\xi) = ||\xi||L_2(\Omega). \eqno(1.4.a)
$$

 Recall that always $ ||\xi ||\Sub \ge \sigma(\xi) = ||\xi||L_2(\Omega),  $ so that

$$
\xi \in \SSub(\Omega) \Leftrightarrow  {\bf E}\xi = 0, \ ||\xi ||\Sub = \sigma(\xi) = ||\xi||L_2(\Omega).  \eqno(1.4b)
$$

 Many examples of strictly subgaussian distributions may be found in the book of  V.V.Buldygin and Yu.V.Kozachenko \cite{Buldygin2},
chapter 1. For instance, arbitrary mean zero Gaussian distributed r.v. is strictly subgaussian, including the case when
this r.v. is equal to zero  a.e.; the symmetric Rademacher's r.v. $ \rho $ with distribution
$ {\bf P}(\rho = 1) = {\bf P}(\rho = -1) = 1/2  $  belongs to the set  $ \SSub(\Omega). $
The random variable $  \eta $ which has an uniform distribution on the symmetrical interval $ (-b,b), \ b = \const \in (0,\infty) $ is
Strictly Subgaussian. \par
 Consider also following the authors \cite{Buldygin2}  the r.v. $ \zeta $  with the following density:

 $$
 f_{\zeta}(x) =   \frac{\alpha+1}{2 \alpha} \  \left(1 - |x|^{\alpha} \right) \ I(|x| \le 1), \
 \alpha = \const \ge 0, \eqno(1.5)
 $$
 where $ I(A) = I(A,x) = 1, \ x \in A; \ I(A) = I(A,x) = 0, \ x \notin A $ is indicator function;
then $ \zeta \in \SSub(\Omega). $ \par

This example is interesting because the kurtosis of the r.v. $ \zeta $ is zero if $ \alpha = \sqrt{10} - 3. $

\vspace{3mm}

 The convenience of these notions is following. Let $  \{ \xi(i ) \}, \ i = 1,2,\ldots,n  $ be (centered) independent
subgaussian r.v. Denote

$$
S(n) = \sum_{i=1}^n \xi(i), \hspace{6mm} \Sigma^2(n) = \sum_{i=1}^n (||\xi(i)||\Sub)^2. \eqno(1.6)
$$

 Then $ ||S(n)||\Sub \le \Sigma(n) $ and following

$$
 \max( {\bf P}(S(n)/\Sigma(n) > x ),  {\bf P}(S(n)/\Sigma(n) < -x )) \le e^{ -x^2/2  }, \ x \ge 0, \eqno(1.7)
$$
 the tail or concentrations inequalities. \par
  If in addition $ \xi(i) $ are identical distributed and $ \beta:= ||\xi(1)||\Sub \in (0,\infty), $ then

 $$
 \sup_n ||S(n)/\sqrt{n}||\Sub = \beta
 $$
  and

$$
\sup_n \max( {\bf P}(S(n)/(\beta \sqrt{n} ) > x ),( {\bf P}(S(n)/(\beta \sqrt{n} ) < - x ) \le e^{ -x^2/2  }, \ x \ge 0, \eqno(1.8)
$$

 If in addition the r.v. $ \xi(i) $  are strictly subgaussian, the estimate (1.8) may be reinforced by lower estimate  used the classical CLT:

$$
\sup_n {\bf P}(S(n)/(\beta \sqrt{n} ) > x ) \ge \lim_{n \to \infty} {\bf P}(S(n)/(\beta \sqrt{n} ) > x ) =
$$

$$
( 2 \pi )^{-1/2} \int_x^{\infty} e^{-y^2/2} \ dy \ge C  \ x^{-1} e^{-x^2/2}, \ x \ge 1. \eqno(1.9)
$$

\vspace{3mm}

{\bf This short report may be considered as a slight addition to the book of V.V.Buldygin and Yu.V.Kozachenko
\cite{Buldygin2}; we give some new examples of Subgaussian and
Strictly Subgaussian random variables, obtain the exponential exact bounds for tails of distribution for
sums of independent Strictly Subgaussian random variables, extend this estimates on the sequence of martingale differences,
 investigate the mixture of ones distributions etc.   }\par

\vspace{3mm}

 Applications of these notions in the non-parametrical statistics may be found in the articles
\cite{Gaivoronsky1}, \cite{Kiefer1}. Another statistical applications is described in \cite{Chen1}, \cite{Ryiabinin1}.
The subgaussian r.v. appears also in the articles \cite{B'enyi1}, \cite{B'enyi2} devoted to the non-linear
Schr\"odinger's equation.  Some applications in the information and coding theory  see  in \cite{Raginsky1}. \par

 \vspace{3mm}

 \section{Mixture of subgaussian random variables}

 \vspace{3mm}

 Let $  Z = \{  z  \}  $ be another set equipped some sigma-algebra and probability measure $ \mu. $ Let also
 $  \xi_z, \ z \in Z $ be a {\it  a family  } of random variables  such that the function

 $$
 z \to {\bf P} (\xi_z \in A), \ A \in B
 $$
is $ \mu \ - $ measurable. By definition, the random variable $  \nu, $ more precisely its distribution is called
{\it mixture  } of individual distributions $ {\bf P} ( \xi_z \in A ) $ relative the (weight) measure  $ \mu, $ if

$$
{\bf P} (\nu \in A) = \int_Z {\bf P} ( \xi_z \in A ) \ \mu(dz).  \eqno(2.1)
$$

Another interpretation-conditional distribution.\par

  If the equality (2.1) there holds, then for all non-negative measurable function $ h: R \to R $

$$
{\bf E} h(\nu) = \int_Z {\bf E} h(\xi_z) \ \mu(dz). \eqno(2.2)
$$

\vspace{3mm}

{\bf Theorem 2.1.}  Let the random variables $  \xi_z $ be Strictly Subgaussian with {\it at the same subgaussian norm}

$$
||\xi_z||\Sub = \sigma = \const \in (0, \infty), \eqno(2.3)
$$
  Then the random variable $ \nu $ is also Strictly Subgaussian with at the same norm. \par

\vspace{3mm}

{\bf Proof.} We derive using the identity (2.2) for the functions $ h(x) = x, x^2, e^{\lambda x} $ correspondingly

$$
{\bf E} \nu = \int_Z {\bf E} \xi_z \ \mu(dz) = 0, \  {\bf E} \nu^2 = \int_Z {\bf E} \xi^2_z \ \mu(dz) = \int_Z  \sigma^2 \ \mu(dz)=
\sigma^2,
$$

$$
{\bf E} e^{\lambda \nu} = \int_Z {\bf E} e^{\lambda \xi_z} \ \mu(dz) \le \int_Z e^{\lambda^2 \sigma^2/2} \ \mu(dz) =  e^{\lambda^2 \sigma^2/2},
$$
Q.E.D.\par

Where the conditions of Theorem 2.1 are not satisfied, the variable $  \nu  $ may be as a Strictly Subgaussian
or no. The correspondent (very spectacular)  example  see in the aforementioned monograph \cite{Buldygin2}, p. 14:

$$
{\bf P} (\nu = 1) = {\bf P}(\nu = -1) = \frac{1-\gamma}{2}, \ {\bf P}(\nu = 0) = \gamma, \ \gamma = \const \in [0,1];
$$
then $  \nu \in \SSub(\Omega) $ iff $  0 \le \gamma \le 2/3 $ or $ \gamma = 1; $ otherwise $ \nu \notin \SSub(\Omega). $\par

 Let now $ \xi_1 $ be a centered Gaussian distributed r.v. with variance $ \sigma^2_1 $ and $  \xi_2 $ be also mean zero Gaussian variable with
other variance $ \sigma_2^2,   $ wherein  $ 0 < \sigma_1^2 < \sigma^2_2 < \infty   $ and $  \mu (\{  1\}) = \mu ( \{ 2 \} )= 1/2. $ It is easy to
verify that the r.v. $ \nu $ is not strictly subgaussian despite both the r.v. $ \xi_1, \ \xi_2  $ are strictly subgaussian.\par

\vspace{3mm}

{\bf Remark 2.1.} If the random variables $ \xi(z) $ in the theorem 2.1  are only subgaussian with variable but uniformly bounded norm

$$
\sigma_+ := \sup_z ||\xi_z||\Sub =  \sup_z \sigma(z) \in (0,\infty), \eqno(2.4)
$$
then the r.v. $ \nu  $ is subgaussian with subgaussian norm less than $ \sigma_+. $\par

 \vspace{3mm}

 \section{Disjoint  subgaussian random variables}

 \vspace{3mm}

 Two r.v. $ \eta_1, \eta_2 $ are called by definition {\it disjoint,} if $ \eta_1 \cdot \eta_2 = 0 $ almost everywhere. The
{\it family } $  \eta_j, \ j = 1,2,\ldots,n; \ n  < \infty $ of r.v.   is named disjoint, if it is pairwise disjoint. \par

 Denote  $ \beta(j) = ||\eta_j||\Sub, \  S(n) = \sum_{j=1}^n \eta_j  $ and suppose $  {\bf E} \eta_j = 0. $\par

\vspace{3mm}

 {\it We intend to obtain in this section the subgaussian and other exponential estimations for the sums of disjoint random variables.} \par

\vspace{3mm}

 The first $  L_p(\Omega) $  estimates for $ S(n) $ was obtained in the famous work of H.P.Rosenthal \cite{Rosenthal1};
the modern results with very interest generalization
see in the articles of  S.V.Astashkin and  F.S.Sukochev \cite{Astashkin1}, \cite{Astashkin2}.\par

\vspace{3mm}

 Let us introduce the following function:

$$
G_n(y_1,y_2, \ldots, y_n) \stackrel{def} =
\inf_{ \mu > 0} \left[ \mu^{-1} \ln \left(  \sum_{j=1}^n e^{\mu y_j} - (n-1)  \right) \right]^{1/2}. \eqno(3.1)
$$

\vspace{3mm}

{\bf Theorem 3.1.} {\it Let $  \eta_j, \ j = 1,2,\ldots,n $ be centered disjoint random variables. Then}

$$
||S(n)||\Sub  \le G_n(\beta_1^2,\beta_2^2, \ldots, \beta_n^2). \eqno(3.2)
$$

\vspace{3mm}

{\bf Proof.} Given:  $  \eta_j = \eta_j \cdot I(A(j)), \ A(j) \cap A(i) = \emptyset, i \ne j, \ A(j) \in B; $

$$
\int_{\Omega} e^{ \lambda \eta_j  } {\bf P}(d \omega) \le e^{\lambda^2 \beta^2_j/2},
$$
therefore

$$
\int_{A(j)}e^{ \lambda \eta_j } {\bf P}(d \omega) \le e^{ \lambda^2 \beta^2_j/2 } - (1 - {\bf P}(A_j)).
$$

 We deduce:

 $$
 {\bf E}e^{\lambda S(n) }= {\bf E} e^{  \lambda \sum_j \eta_j } = {\bf E} e^{  \lambda \sum_j \eta_j  I(A(j))}=
 {\bf E} \prod_j e^{\lambda \eta_j I(A(j)) } =
 $$

$$
\sum_j \int_{A(j)} e^{ \lambda \eta_j } {\bf P}(d \omega) + (1 - \sum_j {\bf P}(A(j))) \le
$$

$$
\sum_j \left[ e^{ \lambda^2 \beta^2_j/2  } -( 1 - {\bf P}(A(j))  ) \right]  + (1 - \sum_j {\bf P}(A(j)))  =
$$

$$
\sum_j e^{ \lambda^2 \beta^2_j/2  } -(n-1) \le e^{ \lambda^2 G^2_n(\beta_1^2, \beta_2^2, \ldots, \beta^2_n) /2 } \eqno(3.3)
$$
on the basis of definition of the function $  G_n(\cdot). $  Theorem 3.1 is proved. \par

 \vspace{3mm}

{\bf  Remark 3.1. } It is interest to note that  our estimation does not depend on the partition

$$
 R = \{ A(j), j = 1,2,\ldots,n; \ \Omega \setminus \cup_j A(j)   \}.
$$

\vspace{3mm}

{\bf  Example 3.1. } Let $ \nu:  \Omega \to R $ be a centered stepwise (simple) r.v. (measurable function):

$$
\nu = \sum_{j=1}^{m} c(j) [ I(A(p(j))) - p(j)  ], \ m = \const \le \infty,  \ c(j) = \const,
$$
and $ \{ A(p(j)) \} $ are pairwise disjoint events. We conclude  using triangle inequality
for the subgaussian norm and the completeness of the space $  \Sub(\Omega) $  in the case when $  m = \infty: $

$$
||\nu||\Sub \le \sum_{j=1}^m |c(j)| Q(p(j)).
$$
 The application of theorem 3.1 gives more exact estimation. \par

\vspace{3mm}

{\it  Another approach to the problem of exponential tail estimates for sums of disjoint random variables.}\par

\vspace{3mm}

 We recall briefly first of all here  for reader conventions some definitions and facts from
the theory of GLS spaces.\par

 Recently, see \cite{Fiorenza1}, \cite{Fiorenza2},\cite{Ivaniec1}, \cite{Ivaniec2}, \cite{Jawerth1},
\cite{Karadzov1}, \cite{Kozachenko1}, \cite{Liflyand1}, \cite{Ostrovsky1}, \cite{Ostrovsky2} etc.
 appear the so-called Grand Lebesgue Spaces (GLS)
 $$
 G(\psi) = G = G(\psi ; B);  \ B = \const \in (1, \infty]
 $$
spaces consisting on all the random variables  (measurable functions) $ f : \Omega \to R  $ with finite norms

$$
||f||G(\psi) \stackrel{def}{=} \sup_{p \in (A;B)} \left[\frac{|f|_p}{\psi(p)} \right], \eqno(3.4)
$$

$$
|f|_p  \stackrel{def}{=} \left[ {\bf E} |f|^p \right]^{1/p}, \ 1 \le p \le \infty.
$$

 Here $ \psi = \psi(p), \ p \in [1,B) $ is some continuous positive on the {\it open} interval $ (1;B) $ function such
that

$$
\inf_{p \in(A;B)} \psi(p) > 0. \eqno(3.5)
$$

We will denote
$$
\supp(\psi) \stackrel{def}{=} [1;B)
$$
 or by abuse of notations $ \supp(\psi) = B.  $ \par

The set of all such a functions with the support $ \supp(\psi) = (1;B) $ will be denoted by  $  \Psi(1;B) = \Psi(B). $  \par

This spaces are rearrangement invariant; and are used, for example, in
the theory of Probability, theory of Partial Differential Equations,
 Functional Analysis, theory of Fourier series,
 Martingales, Mathematical Statistics, theory of Approximation  etc. \par

 Notice that the classical Lebesgue - Riesz spaces $ L_p $  are extremal case of Grand Lebesgue Spaces, see
 \cite{Ostrovsky2}. \par

 Let a function $  \xi:  \Omega \to R  $ be such that

 $$
 \exists B > 1  \Rightarrow  \forall p \in [1, B)  \ |\xi|_p < \infty.
 $$
Then the function $  \psi = \psi_{\xi}(p) $ may be naturally defined by the following way:

$$
\psi_{\xi}(p) := |\xi|_p, \ p \in [1,B). \eqno(3.6)
$$

 The finiteness of the $ G\psi \ - $ norm for some r.v. $  \xi $ allows to obtain the exact exponential tail inequalities for the
distribution  $  \xi; $ for instance,

$$
\sup_{p \ge 1} \left[ \frac{|\xi|_p}{p^m} \right]  < \infty \ \Leftrightarrow \exists C > 0, \ \forall x \ge 0 \ \Rightarrow
{\bf P} (|\xi| > x) \le e^{ - C x^{1/m} }, \ m = \const > 0, \eqno(3.7)
$$
see \cite{Kozachenko1}, \cite{Ostrovsky1}, chapter 1, section 3.\par

\vspace{3mm}

 Let  us return to the the problem of exponential estimations for sums of disjoint variables.\par

\vspace{3mm}

{\bf  Proposition 3.1.  } Let $ \xi, \eta  $ be two disjoint r.v. belonging to some space $ G\psi, \ \psi \in G\Psi $ with
$ \supp \psi = B \in (1, \infty]. $ Then

$$
||\xi + \eta||G\psi \le \left[ (||\xi||G\psi)^B +   ( ||\eta||G\psi)^B   \right]^{1/B}, \eqno(3.8)
$$
where at $  B = \infty $

$$
||\xi + \eta||G\psi \le \max [ ||\xi||G\psi,  ||\eta||G\psi ]. \eqno(3.8a)
$$

{\bf Proof.} Denote for brevity $ a = ||\xi||G\psi, \    b = ||\eta||G\psi. $ Then

$$
|\xi|_p \le a \psi(p), \ |\eta|_p \le b \psi(p); \  |\xi|_p^p \le a^p \psi^p(p), \ |\eta|_p^p \le b^p \psi^p(p).
$$
Since the r.v.  $  \xi, \eta $ are disjoint,

$$
|\xi + \eta|_p^p = |\xi|_p^p + |\eta|_p^p = (a^p + b^p) \cdot \psi^p(p).
$$

 Therefore,

$$
|\xi + \eta|_p \le \psi(p) \cdot (a^p + b^p)^{1/p} \le \psi(p) \cdot (a^B + b^B)^{1/B}.
$$
 It remains to divide on the $ \psi(p) $ and take maximum over $ p; 1\le p < B. $\par

 The generalization on the sum of $  n  $ disjoint variables is clear.\par
%&

\vspace{3mm}

 \section{Martingale case}

 \vspace{3mm}

 Let $  F $ be non-trivial sigma subalgebra of the source sigma-algebra $  B. $  The r.v. $ \eta $ is said to be
{\it conditional subgaussian,}  if there is a {\it non-random } non-negative constant $  \tau  $ for which

$$
\forall \lambda \in R \ \Rightarrow   {\bf E} e^{ \lambda \eta  }/F \le e^{ \lambda^2 \tau^2/2  }. \eqno(4.1)
$$

 The minimal value of the constant $  \tau $ from the inequality (4.1) is called {\it conditional subgaussian norm}
of the r.v. $  \eta $  relative the sigma-algebra $  F, $ write

$$
||\eta||\Sub(F):=  \inf \{ \tau, \ \tau > 0: \
 \forall \lambda \in R \ \Rightarrow   {\bf E} e^{ \lambda \eta  }/F \le e^{ \lambda^2 \tau^2/2  }  \}. \eqno(4.2)
$$

 The set of all r.v. with $ ||\eta||\Sub(F) < \infty   $  relative this norm and ordinary algebraic operations
 forms by definition the complete Banach space $  \Sub(\Omega,F). $ \par

For example, arbitrary centered bounded r.v.  $ \eta = \eta(\omega), \ \omega \in \Omega  $ belongs to any space $ \Sub(\Omega,F).  $\par

Obviously, $  ||\eta||\Sub(\Omega) \le  ||\eta||\Sub(\Omega,F). $\par

 If in the equality the value $ \tau $ may  be selected such that

 $$
\tau = ||\eta||\Sub(F) = \sqrt{ \Var(\eta)/F } = \sqrt{ {\bf E} \eta^2/F  },
 $$
 then as before the r.v. $  \eta  $ may be named Strictly  conditional subgaussian:  $  \eta \in \SSub(\Omega,F). $\par

\vspace{4mm}

 Recall that sequence $ (X(i), F(i)), \ i = 0,1,2,\ldots,n  $ where the $  X(i) $ are random
variables and the $ F(i) $ are sigma-algebras, is a martingale if the following
conditions are satisfied: \\

1. The sequence of sigma-algebras $ \{ F(i) \} $ forms a filtration, i.e. $ F(0) \subset F(1) \subset F(2) \ldots \subset F(n); $ usually, $ F(0) $
is the trivial sigma-algebra $  (\emptyset, \Omega)  $  and $ F(n) $ is sigma-subalgebra  of  source sigma-algebra $ B.$ \\

2. $  X(i) \in L_1(\Omega,P) $ and

$$
X_{i-1} \stackrel{a.e.}{=} {\bf E} X(i)/F(i-1),  \ i = 1,2,\ldots,n-1.
$$

 We suppose in the sequel $ {\bf E} X(i)  = 0 $ and introduce the correspondent sequence of martingale-differences $ \xi(i) $ as follows:

$$
\xi(0) = 0, \ \xi(i) = X(i+1) - X(i), \ i = 1,2,\ldots, n-1.
$$

\vspace{3mm}

 Denote
 $$
  \theta(j) = ||\xi(j)||(\Sub(\Omega, F(j-1)), \ j = 1,,2,\ldots,n, \ \Delta(n) = \sqrt{\sum_{j=1}^n \theta^2(j)  }, \eqno(4.3)
 $$
if there exists.\par

\vspace{3mm}

{\bf Theorem 4.1.}

$$
||X(n)||\Sub(\Omega) \le \Delta(n); \eqno(4.4)
$$

$$
\max( {\bf P}(X(n)/\Delta(n)> x ),  {\bf P}(X(n)/\Delta(n) < - x ) \le e^{-x^2/2}, \ x \ge 0. \eqno(4.4a)
$$

\vspace{3mm}

{\bf Proof }  is alike to the one in the article of K.Azuma \cite{Azuma1}; see also \cite{Raginsky1}. Namely, let $  \lambda = \const \in R; $ then

$$
{\bf E} e^{\lambda X(n)} = {\bf E} e^{ \lambda \sum_{j=1}^{n-1} \xi(j) } =
{\bf E} \left[ {\bf E} \left[ e^{ \lambda \sum_{j=1}^{n-1} \lambda \xi(j) } \right] /F(n-2)  \right] =
$$

$$
{\bf E} \left[  \left[ e^{ \lambda \sum_{j=1}^{n-2} \lambda \xi(j) } \right]  \ \cdot {\bf E} e^{ \lambda \xi(n-1) }/F(n-2)  \right] \le
$$

$$
{\bf E} \left[  \left[ e^{ \lambda \sum_{j=1}^{n-2} \lambda \xi(j) } \right]  \ \cdot  e^{ \lambda^2 \theta^2(n)/2 } \right] \le \ldots \le
$$

$$
e^{ \lambda^2 \sum_{j=1}^n \theta^2(j)/2  } = e^{\lambda^2 \Delta^2(n)/2 }. \eqno(4.5)
$$
Thus,

$$
||X(n)||\Sub(\Omega) \le \Delta(n). \eqno(4.6)
$$
 This completes the proof of theorem 4.1.\par

 Note that this estimate is exponential exact if for instance $ \{ \xi(j) \} $ are independent and strictly subgaussian.\par

\vspace{3mm}

 \section{Another examples of subgaussian random variables}

 \vspace{3mm}

{\bf A. Symmetrized beta distribution.}\par

\vspace{3mm}

 Let us consider a symmetrical r.v. $  \xi = \xi_{\alpha, \beta} $ with the density

$$
f(x) = f_{\alpha,\beta}(x) = 0.5 \frac{ |x|^{\alpha - 1} \ ( 1 - |x|)^{\beta - 1}}{ B(\alpha, \beta)} \ I(|x| < 1), \ \alpha, \beta = \const > 0,
\eqno(5.1)
$$

where as usually $  B(\alpha, \beta) $ denotes the beta-function. \par

\vspace{3mm}

{\bf Theorem 5.1.}  {\it If }

$$
B(\alpha,\beta) \le 1,  \eqno(5.2)
$$
{\it then the r.v. $ \xi = \xi_{\alpha, \beta}  $ is strictly subgaussian. } \par

\vspace{3mm}

Note that the condition (5.2) of theorem 5.1 is satisfied if for instance $ \alpha \ge 1, \ \beta \ge 1.  $ \par

\vspace{3mm}

{\bf Proof.}  We have:

$$
{\bf E} \xi^{2k + 1} = 0, \ {\bf E} \xi^{2k} = \frac{B(2k + \alpha, \beta)}{B(\alpha, \beta )}, \ k = 0,1, \ldots,
$$
so that

$$
\sigma^2 := \Var(\xi) = \frac{B(2 + \alpha, \beta)}{B(\alpha, \beta )}= \frac{\alpha(\alpha + 1)}{(\alpha + \beta)( \alpha + \beta + 1) };
$$

$$
{\bf E} e^{\lambda \xi} = \sum_{k=0}^{\infty} \frac{\lambda^{2k}}{(2k)!} \cdot {\bf E} \xi^{2k} = \sum_{k=0}^{\infty}
\frac{B(2k + \alpha, \beta)}{B(\alpha, \beta )} \cdot \frac{\lambda^{2k}}{(2k)!};
$$

$$
e^{\lambda^2 \sigma^2/2} = \sum_{k=0}^{\infty} \frac{\lambda^{2k} \sigma^{2k}}{2^k \ k!}.
$$

 It is sufficient to prove that

$$
\frac{B(2k + \alpha, \beta)}{B(\alpha, \beta )} \cdot \frac{1}{(2k)!} \le \frac{\sigma^{2k}}{2^k \ k!}, \ k = 1,2,\ldots; \eqno(5.3)
$$
see \cite{Buldygin2}, p.8; the case $  k = 0 $ is trivial.\par

 The inequality (5.3) is equivalent to the following:

$$
\frac{B(2k + \alpha, \beta)}{B(\alpha, \beta )} \cdot \frac{1}{(2k)!}  \le  \frac{B^{2k}(\alpha + 2, \beta)}{B^{2k}(\alpha, \beta) \cdot   2^k \cdot k!}.
\eqno(5.4)
$$

 We denote

$$
\theta(k) = \left[\frac{B(2k + \alpha, \beta)}{B(\alpha, \beta )} \cdot \frac{1}{(2k)!} \right] :
\left[ \frac{B^{2k}(\alpha + 2, \beta)}{B^{2k}(\alpha, \beta) \cdot   2^k \cdot k!} \right], \ k = 0,1,\ldots;
\eqno(5.5)
$$
then $ \theta(0) = 1. $ Further, the expression for the fraction $ \theta(k+1)/\theta(k), \ k = 0,1,2, \ldots  $ has a form:

$$
\theta(k+1)/\theta(k) = \frac{\Gamma^2(\alpha) \Gamma^2(\beta)}{\Gamma^2(\alpha + \beta)} \cdot \frac{k+1}{k(2k+1)} \times
$$

$$
\frac{(2k + \alpha)(2k+ \alpha + 1)}{(2k + \alpha + \beta + 1)(2k + \alpha + \beta)} < B^2(\alpha,\beta) < 1;
$$
therefore $  \theta(k) \le 1, \ k = 0,1,2, \ldots, $ Q.E.D.\par

 \vspace{3mm}

{\bf B. Symmetrized gamma distribution.}\par

\vspace{3mm}

 We consider here the r.v. $ \gamma = \gamma_{\alpha, \beta} $ with the (symmetrical) density of distribution

$$
g_{\alpha, \beta}(x) = 0.5 \ \frac{\beta \ |x|^{\alpha} \ e^{- |x|^{\beta}} }{ \Gamma( (\alpha + 1)/\beta )}, \ x \in R, \ \alpha = \const > -1,
\beta = \const > 0. \eqno(5.6)
$$
 We have:   $   {\bf E} \gamma^{2k+1} = 0, \ k = 0, 1, 2,\ldots;   $

$$
{\bf E} \gamma^{2k} = \frac{\Gamma( (\alpha + 2k +1)/\beta  )}{\Gamma( (\alpha + 1)/\beta)}. \eqno(5.7)
$$
 In particular,

$$
\Var(\gamma) = \frac{\Gamma( (\alpha + 3)/\beta  )}{\Gamma( (\alpha + 1)/\beta)}. \eqno(5.7a)
$$

 We suppose in the sequel $  \beta > 2; $ obviously, when $  \beta < 2, $ the r.v. $ \gamma $ is even not subgaussian.\par
The r.v.   $ \gamma  $ is subgaussian iff $ \beta \ge 2. $

\vspace{3mm}

 In order  to formulate the next result, we need to introduce some preliminary notations.

 $$
 \theta = \max(\alpha - \beta + 1,1),  \hspace{6mm} k_0 = \max(1, (\beta - \alpha - 1)/2),
 $$

 $$
 G_1 = \sqrt{\pi}
 \cdot \Gamma \left( \frac{\alpha+1}{\beta}  \right) \cdot e^{- 3/8 } \cdot e^{ -\theta/(2 \beta)  },
 $$

$$
G_2 = 0.25 \ e  \cdot \frac{\Gamma ( ( \alpha + 3  )/\beta )}{ \Gamma((\alpha + 1)/\beta)} \cdot
\left[ 2 \beta^{-1} (1 + \theta/(2 k_0) \right]^{2/\beta},
$$

$$
G = \left[( \max(1, G_1)  ) \cdot G_2 \right]^{\beta/(\beta - 2) },
$$

$$
\zeta(k)= \zeta_{\alpha,\beta}(k) := \Gamma( ( \alpha + 2k + 1 )/\beta ) \cdot
 \frac{2^k \ k!}{(2k)!} \cdot  \frac{\Gamma^{ k-1 }( ( \alpha +1)/\beta )}{\Gamma( (\alpha + 3)/\beta) }.
$$

Recall that we consider the values $ \beta $ greatest than 2. \par

\vspace{3mm}

 {\bf  Theorem 5.2.} {\it If the following conditions }

$$
\forall k < G \ \Rightarrow \zeta_{\alpha,\beta}(k) \le 1 \eqno(5.9)
$$
 {\it there holds, then the r.v. $ \gamma = \gamma_{\alpha,\beta} $ is strictly subgaussian.} \par

\vspace{3mm}

 {\bf Proof.}\par

It is sufficient to prove as before the following inequality for all the values $ k = 0,1,2, \ldots $

$$
\zeta(k)= \zeta_{\alpha,\beta}(k) := \Gamma( ( \alpha + 2k + 1 )/\beta ) \cdot
 \frac{2^k \ k!}{(2k)!} \cdot  \frac{\Gamma^{ k-1 }( ( \alpha +1)/\beta )}{\Gamma( (\alpha + 3)/\beta) } \le 1. \eqno(5.10)
$$

 Since $  \beta > 2, \ \lim_{k \to \infty} \zeta(k) = 0. $  We  conclude omitting complicated computation more exact
estimate  using  Stirling's formula:

$$
k \ge G  \ \Rightarrow  \zeta_{\alpha, \beta}(k) \le 1.
$$
 We derive taking into account the condition (5.9) that (5.10) is true for all the integer values $ k = 1,2,\ldots,  $ Q.E.D.\par

\vspace{3mm}

{\bf Remark 5.1. Example.} Let $  \alpha = \beta \to \infty;  $ then condition (5.9) is satisfied; moreover,

$$
\lim_{\beta \to \infty} \sup_{k \ge 1} \zeta_{\beta, \beta}(k)  = 0.
$$
 Therefore, there exists a value  $ Z   $ such that for all the values $  \alpha = \beta > Z $ the r.v. $  \gamma_{\beta, \beta}  $
is strictly subgaussian.
%&

\vspace{3mm}

{\bf Remark 5.2.} Recall that the condition $  \gamma(1) \le 3 $ or in detail the inequality

$$
\Gamma((\alpha + 5)/\beta) \Gamma( ( \alpha + 1  )/\beta ) \le 3 \Gamma^2( (\alpha + 3   )/\beta ) \eqno(5.11)
$$
is necessary for strictly subgausianness.  Denote $ \epsilon = 1/\beta $ and assume $ \alpha = \const, \
 \epsilon \to 0+ \ \Leftrightarrow \beta \to \infty. $
 Taking into account the behavior of the Gamma function near to the value $ 0+: $

$$
\Gamma(\epsilon) \sim 1/\epsilon, \ \epsilon \to 0+,
$$
we deduce from (5.11) as $  \beta \to \infty $

$$
\frac{1}{\alpha + 5} \cdot  \frac{1}{\alpha+1} \le  \frac{3}{( \alpha + 3 )^2},
$$
which is not true as $ \alpha \to -1 + 0. $\par

 Moreover, the solution of the last inequality subject to the limitation $ \alpha > -1 $ has a form
 $ \alpha \ge 3(\sqrt{3} - 1) \approx 2.19.   $ \par
  Thus, for all the values $  \alpha $ from the interval $ -1 < \alpha < 3(\sqrt{3} - 1)  $ there exists a positive the value $ \beta_0 $
such that for all the values $ \beta > \beta_0 $  the random value $  \gamma = \gamma_{\alpha,\beta} $ is subgaussian but not strictly subgaussian.

\vspace{4mm}

Another (but more simple) {\bf example } of strictly subgaussian r.v. with {\it unbounded } support may be constructed by means of section 2.
Namely, let the r.v. $ \xi_1 $ has an uniform distribution on the symmetric interval  $ [-b,b], \ b = \const > 0 $ and let $  \xi_2 $ has a normal
(Gaussian) mean zero distribution  with at the same variation $ \Var \xi_2 = b^2/3. $ The arbitrary non - trivial mixture of these distributions
has an unbounded support, is not Gaussian and is strictly subgaussian by theorem 2.1. \par

 \vspace{3mm}

\section{ Acknowledgement.}\par

 Authors are very grateful to prof. S.V.Astashkin  and L.Maligranda for sending Your remarkable articles and comments. \\

%\vspace{3mm}

\end{document}